\documentclass[12pt]{amsart}
\allowdisplaybreaks
\usepackage{amssymb}
\usepackage{graphicx}
\usepackage{fullpage}
\usepackage[pdfencoding=auto,colorlinks,linkcolor=,citecolor=]{hyperref}
\usepackage{cite}
\usepackage{dsfont}
\usepackage{eucal}
\usepackage[all,cmtip]{xy}

\newcommand{\Aut}{\mathsf{Aut}}
\newcommand{\End}{\mathsf{End}}
\newcommand{\Ext}{\mathsf{Ext}}
\renewcommand{\ge}{\geqslant}

\newcommand{\im}{\mathsf{Im}}
\newcommand{\Ker}{\mathsf{Ker}}
\renewcommand{\le}{\leqslant}
\newcommand{\mmod}{\mathsf{mod}}
\newcommand{\op}{{\mathsf{op}}}
\newcommand{\Rad}{\mathsf{Rad}}

\newcommand{\bF}{\mathbb{F}}
\newcommand{\sfM}{\mathsf{M}}

\theoremstyle{plain}
\numberwithin{equation}{section}
\newtheorem{lemma}[equation]{Lemma}
\newtheorem{theorem}[equation]{Theorem}

\theoremstyle{definition}

\newtheorem{example}[equation]{Example}

\theoremstyle{remark} 
\newtheorem{remark}[equation]{Remark} 
 
\newtheorem{question}[equation]{Question}

\author{David J. Benson} 
\address{Institute of Mathematics \\ 
University of Aberdeen \\ 
Aberdeen AB24 3UE \\ 
United Kingdom}

\email{d.j.benson@abdn.ac.uk}

\title{Centralisers of finite groups in \\ locally finite simple groups}

\keywords{Locally finite group, centraliser, representation theory}

\subjclass{Primary: 20F50. Secondary: 20C20, 16G20.}

\begin{document}
\begin{abstract}
We answer in the negative a question of Hartley about representations
of finite groups, by constructing examples of finite simple groups
with arbitrarily large representations whose endomorphism ring
consists of just the scalars. We show as a consequence that there are finite
simple groups of automorphisms of the locally finite simple group $SL(\infty,\bF_q)$
with trivial centraliser. The smallest of our examples is $A_6$ with $q=9$.
\end{abstract}

\maketitle

\section{Introduction}

A group is said to be \emph{locally finite} if every finitely  
generated subgroup is finite. A group is said to be \emph{linear} if  
if it has a faithful representation by matrices over a field.  
 
In the introduction to~\cite{Hartley:1990a}, and Section~3
of~\cite{Hartley:1992a}, Brian Hartley mentioned the following
problem. 

\begin{question}\label{qu:Hartley}
  Is it the case that in a non-linear locally finite simple group, the
  centraliser of every finite subgroup is infinite?
\end{question}

Papers that have investigated Hartley's Question~\ref{qu:Hartley} include 
Brescia and Russo~\cite{Brescia/Russo:2019a}, 
Ersoy and Kuzucuo\u{g}lu~\cite{Ersoy/Kuzucuoglu:2012a}, 
Hartley and Kuzucou\u{g}lu~\cite{Hartley/Kuzucuoglu:1991a}, 
Kuzucuo\u{g}lu~\cite{Kuzucuoglu:1994a,Kuzucuoglu:1997a, 
Kuzucuoglu:2013a,Kuzucuoglu:2017a}. In particular, it is shown 
in~\cite{Hartley/Kuzucuoglu:1991a} that if the finite subgroup is 
cyclic then the answer is ``yes.'' 

A closely related problem is Problem~3.15 of Hartley~\cite{Hartley:1995a}, which asks
the following (with notation slightly altered to align with ours). 

\begin{question}\label{qu:Hartley3.15}
Let $H$ be a finite simple group, or perhaps a finite group such that
$O_p(H)=1$. Suppose that $H\le PSL(n,\bF_q)=K$, where $q$ is a power of
$p$ and $C_K(H)=1$. Does it follow that $n$ is bounded in terms of
$|H|$? We are particularly interested in the case when $H$ is also a
projective special linear group over $\bF_q$. What about the algebraic
group situation, that is, say, $H\cong PSL_m(\bar\bF_p)$, $K\cong
PSL_n(\bar\bF_p)$, with $C_K(H)=1$? Does it follow that $n$ is bounded
in terms of $m$?
\end{question}

An example of a non-linear simple locally finite group is
$SL(\infty,\bF_q)$, the union of the groups $SL(n,\bF_q)$ under the
natural inclusions $SL(n,\bF_q)\to SL(n+1,\bF_q)$ 
fixing the new basis vector. 
It is simple because (except for $n=2$, $q=2$ or $3$) 
every proper normal subgroup of $SL(n,\bF_q)$
consists of scalar multiples of the identity, and are no longer
scalar in $SL(n+1,\bF_q)$ unless trivial. In this group, the
centraliser of a finite subgroup is obviously infinite.

We do not succeed in answering Question~\ref{qu:Hartley},
but we answer a closely related question with the following theorem.

\begin{theorem}\label{th:main}
Let $K_q$ be the locally finite group $SL(\infty,\bF_q)$.
For every prime $p$, there exists a power $q$ of $p$ 
and a finite simple group of automorphisms
$H\le \Aut(K_q)$ such that the centraliser of $H$ in $K_q$ is
trivial.
\end{theorem}

We prove Theorem~\ref{th:main} as a consequence of the following theorem, together
with some information about extensions between simple modules. 
This theorem gives a negative answer to Question~\ref{qu:Hartley3.15},
as is shown in Section~\ref{se:3Ext-eg}. We restrict our attention to
finite groups and linear algebraic groups in order to address
Hartley's question, but it should be clear from the methods that other
situations are also covered.

\begin{theorem}\label{th:Mm}
Let $k$ be a field of characteristic $p$.
Let $H$ be a finite group, or a linear algebraic group over $k$.
Suppose that 
\begin{enumerate}
\item[\rm (i)] $H$
has non-isomorphic simple modules $S \not\cong T$
over $k$ such that
$\End_{H}(S)\cong\End_{H}(T)\cong k$, and
$\dim_k\Ext^1_H(S,T)\ge 3$, \textbf{or}
\item[\rm (ii)] $H$ has non-isomorphic simple modules $R$, $S$ and $T$ with
$\End_{H}(R)\cong\End_{H}(S)\cong\End_{H}(T)\cong k$, 
$\dim_k\Ext^1_{H}(R,T)\ge 1$ and $\dim_k\Ext^1_{H}(S,T)\ge 2$, \textbf{or}
\item[\rm (iii)] $H$ has non-isomorphic simple modules $R$,
  $S_1,\dots, S_r$, and $T_1,\dots,T_r$ with endomorphism rings
  isomorphic to $k$, and
such that the groups
$\Ext^1_{H}(R,T_1)$, $\Ext^1_{H}(S_i,T_i)$,
$\Ext^1_{H}(S_{i+1},T_i)$ and $\Ext^1_{H}(S_1,T_r)$ are all non-trivial.
\end{enumerate}
Then there exists a sequence of indecomposable 
$H$-modules $M_m$ over $k$ and embeddings $M_m\to
M_{m+1}$ such that
\begin{enumerate}
\item
The radical (i.e., the intersection of the maximal submodules) 
$\Rad(M_m)$ is isomorphic to a direct sum of $m$ copies of $T$.
\item
$M_m/\Rad(M_m)$ is isomorphic to a direct sum of $m$ copies of $S$ 
in case {\rm (i)} and a direct sum of one copy of $R$ and $m-1$ copies of $S$
in case {\rm (ii)}.
\item
The endomorphism ring $\End_{H}(M_m)$ consists just of scalar 
multiples of the identity map.
\end{enumerate}
\end{theorem}

The three cases of Theorem~\ref{th:Mm} are proved in
Sections~\ref{se:3Ext}, \ref{se:2Ext} and~\ref{se:1Ext}, and examples
illustrating these three cases are given in Sections~\ref{se:3Ext-eg},
\ref{se:2Ext-eg} and~\ref{se:1Ext-eg}. Theorem~\ref{th:main} is proved
in Section~\ref{se:constr}.\medskip

\noindent
{\bf Acknowledgements.} I would like to thank Alexandre Zalesskii for
introducing me to these problems at breakfast at the M{\o}ller
Institute in Cambridge. This was
while attending the post Covid-19 resumption of the programme 
\emph{`Groups, representations and applications: new
perspectives’} in 2022, at the Isaac Newton Institute,
supported by EPSRC grant EP/R014604/1. My thanks also go to the Isaac
Newton Institute for their support and hospitality.

\section{\texorpdfstring{Three dimensional $\Ext^1$ groups}
{Three dimensional Ext¹ groups}}\label{se:3Ext}

In this section, we deal with case (i) of Theorem~\ref{th:Mm}.
Let $k$ be a field, and 
let $Q$ be the quiver with two vertices and three arrows given by the following diagram.
\[ \xymatrix{ \bullet \ar@<1ex>[r]\ar[r]\ar@<-1ex>[r]& \bullet } \]
A representation of $Q$ over $k$ consists of two
$k$-vector spaces $V$ and $W$ together with three $k$-linear maps
$ X,Y,Z \colon V \to W$. 
We are interested in properties of a particular 
family of representations $\sfM_m$ of $kQ$ for $m\ge 1$ and embeddings 
$\sfM_m\to \sfM_{m+1}$.

For the representation $\sfM_m$, the spaces $V$ and $W$ have dimension $m$, 
and bases $v_1,\dots,v_m$ and $w_1,\dots,w_m$. The actions of $X$, $Y$
and $Z$ on $v_i$ given as follows.
\begin{equation}\label{eq:action}
\begin{array}{|c|ccc|} \cline{2-4}
\multicolumn{1}{c|}{}&X&Y&Z \\ \hline
i=1&w_1&0&0 \\
2\le i\le m&0&w_i&w_{i-1} \\ \hline
\end{array} 
\end{equation}
This can be pictured as follows.
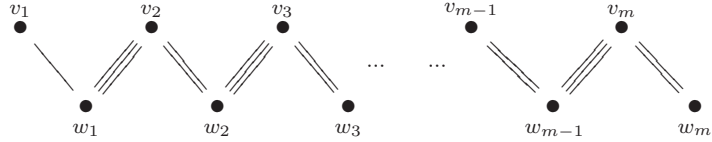
\[ \xymatrix@=3mm{\overset{v_1}{\bullet}\ar@{-}[dr]&&
\overset{v_2}{\bullet}\ar@{=}[dr]\ar@3{-}[dl]&&
\overset{v_3}{\bullet}\ar@{=}[dr]\ar@3{-}[dl]&&\ar@{}[dl]|\cdots
&\overset{v_{m-1}}{\bullet}\ar@{=}[dr]\ar@{}[dl]|\cdots&&
\overset{v_m}{\bullet}\ar@{=}[dr]\ar@3{-}[dl]\\
&\underset{w_1}{\bullet}&&\underset{w_2}{\bullet}
&&\underset{w_3}{\bullet}&&&\underset{\!\!\!\! w_{m-1}\!\!\!\!\!}{\bullet}
&&\underset{w_m}{\bullet}} \]
Here, the actions of $X$, $Y$ and $Z$ are described by the
single, double and triple downward edges respectively.
The embedding $\sfM_m\to \sfM_{m+1}$ takes the basis vectors
$v_i$ and $w_i$ of $\sfM_m$ to the basis elements of $\sfM_{m+1}$ 
with the same names.

\begin{theorem}\label{th:End3}
The endomorphism ring of the $kQ$-module $\sfM_m$ is equal to $k$, 
acting as scalar multiples of the identity map.
\end{theorem}
\begin{proof}
We have direct sum decompositions $V=\Ker(Y) \oplus \Ker(X)$ and
$W=\im(X) \oplus \im(Y)$. Here, $\Ker(Y)$ is spanned by $v_1$ and
$\Ker(X)$ is spanned by $v_2,\dots,v_m$. Similarly,
$\im(X)$ is spanned by $w_1$ and $\im(Y)$ is spanned by
$w_2,\dots,w_m$.

If $\alpha$ is an endomorphism of $\sfM_m$ then $X\alpha(v_1)=\alpha(Xv_1)=\alpha(w_1)$, 
so $\alpha(w_1)$ is in the image of $X$. It is therefore a multiple of $w_1$.
Subtracting off a multiple of the identity map, we can suppose that $\alpha(w_1)=0$. 
Then $\alpha(v_1)$ is in $\Ker(X)\cap\Ker(Y)$ and is hence also zero.

Suppose by induction on $j$ that we have shown that we have 
$\alpha(v_i)=0$ and $\alpha(w_i)=0$ for $1\le i < j$. Then 
$Y\alpha(v_j)=\alpha(Yv_j)=\alpha(w_{j-1})=0$ and
$X\alpha(v_j)=\alpha(Xv_j)=0$, so $\alpha(v_j)\in \Ker(Y)\cap
\Ker(X)$ and hence $\alpha(v_j)=0$. Then
$\alpha(w_j)=\alpha(Yv_j)=Y\alpha(v_j)=0$. 
It follows by induction that for $1\le j \le m$ we have
$\alpha(v_j)=0$ and $\alpha(w_j)=0$, and hence $\alpha=0$.
\end{proof}

\begin{proof}[Proof of Theorem~\ref{th:Mm} in case {\rm (i)}]
Suppose that $k$ is a field of characteristic $p$ and $H$ is a finite
group, or a reductive algebraic group,
with non-isomorphic simple (left) modules $S$ and $T$ 
satisfying $\End_{H}(S)=\End_{H}(T)=k$ and 
$\dim_k\Ext^1_{H}(S,T)\ge 3$. 
Then there is an $H$-module $\Delta$
satisfying $\Rad(\Delta)\cong T \oplus T \oplus T$ and $\Delta/\Rad(\Delta)\cong S$.
Setting $B=\Delta\oplus T$, 
the algebra $\End_{H}(B)^\op$ is isomorphic to
the quiver algebra $kQ$ described above. So $B$ has a right 
action of $kQ$ commuting with the left action of $H$. As a right
$kQ$-module, $B$ is a projective generator for $\mmod$-$kQ$.
Thus $B$ is a $H$-$kQ$-bimodule with
the property that the functor 
$B\otimes_{kQ} -\colon kQ$-$\mmod \to H$-$\mmod$
is fully faithful and exact. It
sends the simple $kQ$-modules to $S$ and $T$ and the three arrows in $J(kQ)$ to three
linearly independent elements of $\Ext^1_{H}(S,T)$. 
So $M_m=B \otimes_{kQ} \sfM_m$ 
is a module of the following form:
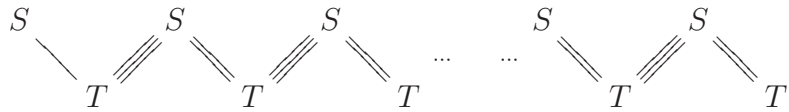
\[ \xymatrix@=5mm{S\ar@{-}[dr]&&
S\ar@2{-}[dr]\ar@3{-}[dl]&&
S\ar@2{-}[dr]\ar@3{-}[dl]&&\ar@{}[dl]|\cdots
&S\ar@2{-}[dr]\ar@{}[dl]|\cdots&&
S\ar@3{-}[dl]\ar@2{-}[dr] \\
&T&&T&&T&&&T&&T} \]
where the three types of lines represent the three  
linearly independent extension classes chosen in the construction of $U$.
By Theorem~\ref{th:End3}, we have $\End_{H}(M_m)\cong k$
and $M_m$ embeds in $M_{m+1}$.
Setting $d=\dim_k S + \dim_k T$, we have $\dim_k M_m=md$.
This completes the proof of Theorem~\ref{th:Mm}.
\end{proof}

\section{\texorpdfstring{Examples involving three dimensional $\Ext^1$ groups}
{Examples involving three dimensional Ext¹ groups}}\label{se:3Ext-eg}

In this section, we indicate that there are many examples of
absolutely irreducible
modules $S$ and $T$ over reductive algebraic groups 
and over finite simple groups $H$, such that
$\Ext^1_{kH}(S,T)$ is at least three dimensional, and sometimes much
larger. The first such examples were discovered by Scott, and already,
feeding these into Case (i) of Theorem~\ref{th:Mm} gives negative
answers to both parts of Question~\ref{qu:Hartley3.15}. 

\begin{theorem}\label{th:Scott}
Suppose that $n\ge 6$. Then for all large enough primes $p$, 
there exist irreducible $PSL(n,\bar\bF_p)$-modules $S$ and $T$ such
that 
\[ \dim_{\bar\bF_p}\Ext^1_{PSL(n,\bar\bF_p)}(S,T) \ge 4. \]
and all large enough powers $q$ of $p$, the restrictions of $S$ and
$T$ are absolutely irreducible 
$\bF_qPSL(n,\bF_q)$-modules such that 
\[ \dim_{\bF_q}\Ext^1_{\bF_qPSL(n,\bF_q)}(S,T)\ge 4. \]
\end{theorem}
\begin{proof}
This is proved in Section~2 of Scott~\cite{Scott:2003a}, by computing
coefficients of Kazhdan--Lusztig polynomials. The
stabilisation of the dimensions of $\Ext^1$ for the finite groups to the dimension
for the algebraic group comes from Theorem~2.8 of
Andersen~\cite{Andersen:1987a}, see Remark~(iii) at the end of that section.
\end{proof}

\begin{remark}
It is shown in Section~4 of L\"ubeck~\cite{Luebeck:2020a} 
that for all large 
enough primes, taking $S$ to be the trivial 
module $\bar\bF_p$, there are simple
modules $T$ for $PSL(n,\bar\bF_p)$ with 
\[ H^1(PSL(n,\bar\bF_p),T)\cong\Ext^1_{PSL(n,\bar\bF_p)}(\bar\bF_p,T) \]
having the following dimensions,
\[ \renewcommand{\arraystretch}{1.4}
\begin{array}{|c|cccc|} \hline
n&6&7&8&9 \\ \hline
\dim_{\bF_q}H^1(PSL(n,\bar\bF_p),T)&3&16&469&36672 \\ \hline
\end{array}\medskip \]
and then for large enough powers $q$ of $p$, the restriction of these $T$ 
are irreducible modules for $\bF_qPSL(n,\bF_q)$ with these dimensions for
$H^1(PSL(n,\bF_q),T)$. 
It is not known whether these dimensions are unbounded
for larger values of $n$, though that seems very likely. 
In the same paper, L\"ubeck proves similar results for 
the other classical groups of types $B_n$, $C_n$ and $D_n$, 
as well as groups of types $F_4$ and $E_6$.
These groups could therefore also be used equally well 
as examples in Theorem~\ref{th:main}. There are also cross
characteristic examples of large dimensional $\Ext^1$ 
for simple modules over groups of Lie type, but the Lie rank in these
cases needs to be a lot larger.

The problem with the theorems above is that the phrase ``large
enough'' is difficult to quantify. 
Some more explicit examples of three dimensional $H^1(H,T)$ for
$T$ an 
absolutely irreducible module for a finite simple group $H$ can be found
in Bray and Wilson~\cite{Bray/Wilson:2008a}, 
for the group $PSU(4,\bF_3)$ over $\bF_3$ with $\dim_{\bF_3}T=19$
(easily verified with \textsc{Magma}~\cite{Bosma/Cannon/Playoust:1997a}),
and for the group
${^2\hspace{-0.2ex}}E_6(\bF_2)$ over $\bF_2$ with $\dim_{\bF_2}
T=1702$
(not so easily verified).
\end{remark}

\section{\texorpdfstring{Two dimensional $\Ext^1$ groups}
{Two dimensional Ext¹ groups}}\label{se:2Ext}

In this section, we deal with case (ii) of Theorem~\ref{th:Mm}. The
proof is very much like that of case (i).
We begin with the quiver $Q$ with three vertices and three arrows as
follows.
\[ \xymatrix{ \bullet \ar@<0.5ex>[r]\ar@<-0.5ex>[r] &
    \bullet & \bullet\ar[l]} \]
A representation of $kQ$ over a field $k$ consists of three vector
spaces $U$, $V$, and $W$ together with three $k$-linear maps\vspace{-3mm}
\[ \xymatrix{V \ar@<0.5ex>[r]^X\ar@<-0.5ex>[r]_Y &
W & U\ar[l]_Z} \]
For the $kQ$-module $\sfM_m$, $U$ is spanned by $v_1$, $V$ is
spanned by $v_2,\dots,v_m$, and $W$ is spanned by $w_1,\dots,w_m$.
The actions of $X$, $Y$ and $Z$ are exactly the same as in 
Table~\eqref{eq:action}.

\begin{theorem}\label{th:End2}
The endomorphism ring of the $kQ$-module $\sfM_m$ is equal to $k$, 
acting as scalar multiples of the identity map.
\end{theorem}
\begin{proof}
The proof is the same as for Theorem~\ref{th:End3}.
\end{proof}

\begin{proof}[Proof of Theorem~\ref{th:Mm} in case {\rm (ii)}]
Suppose that $k$ is a field of characteristic $p$, and let 
$H$ be a finite group with non-isomorphic simple modules $R$, $S$ and $T$
with $\End_{kH}(R)=\End_{kH}(S)=\End_{kH}(T)=k$, $\dim_k\Ext^1_{kH}(R,T)\ge 1$,
and $\dim_k\Ext^1_{kH}(S,T)\ge 2$.
Then there are $kH$-modules $\Delta_1$ and $\Delta_2$ with
\[ \Rad(\Delta_1)\cong T,\qquad \Delta_1/\Rad(\Delta_1)\cong R,\qquad
\Rad(\Delta_2)\cong T\oplus T, \qquad\Delta_2/\Rad(\Delta_2)\cong S. \]
Setting $B=\Delta_1\oplus \Delta_2\oplus T$,  the algebra
$\End_{kH}(B)^\op$ is isomorphic to $kQ$. 
This makes $B$ a
$kH$-$kQ$-bimodule with the property that as a right $kQ$-module it
is a projective generator for $\mmod$-$kQ$.
The functor $B \otimes_{kQ} -$ is fully faithful and exact. It sends
the three simple $kQ$-modules to $R$, $S$ and $T$. So this time,
$M_m=B\otimes_{kQ}\sfM_m$ is a module of the following form\vspace{-2mm}
\[ \xymatrix@=5mm{R\ar@{-}[dr]&&
S\ar@2{-}[dr]\ar@3{-}[dl]&&
S\ar@2{-}[dr]\ar@3{-}[dl]&&\ar@{}[dl]|\cdots
&S\ar@2{-}[dr]\ar@{}[dl]|\cdots&&
S\ar@3{-}[dl]\ar@2{-}[dr] \\
&T&&T&&T&&&T&&T} \]
whose endomorphism ring is $k$, acting by scalar multiples of the identity.
\end{proof}

\section{\texorpdfstring{Examples involving two dimensional $\Ext^1$ groups}
{Examples involving two dimensional Ext¹ groups}}\label{se:2Ext-eg}

In this section, we give examples of finite simple groups $H$ and
simple modules $R$, $S$ and $T$ satisfying the hypotheses of
Case~(ii) of Theorem~\ref{th:Mm}.

\begin{example}\label{eg:p2}
Let $H$ be the group $PSL(2,q)$ with $q=p^2$, $p$ odd, 
and let $k=\bF_q$. Then by 
Corollary~4.5 of Andersen, J{\o}rgensen and Landrock,
there are exactly two simple modules $S$ and $T$,
of dimensions $(\frac{p-1}{2})^2$ and $(\frac{p+1}{2})^2$, such
that $\Ext^1_{kH}(S,T)$ and $\Ext^1_{kH}(T,S)$ are two dimensional.
All the other $\Ext^1$ groups between simple $kH$-modules 
are zero or one dimensional. 
In fact, the entire quiver with relations in this example may be found  
in Koshita~\cite{Koshita:1998a}.  
Since $S$ and $T$ are in the principal
block, and are not the only simples in the principal block, we can
choose a simple module $R$ not isomorphic to $S$ or $T$, so that
$\Ext^1_{kH}(R,T)\ge 1$. For example, if $p=3$ then $PSL(2,9)\cong
A_6$ has simples of dimensions $3$, $1$ and $4$ for $R$, $S$ and $T$.
\end{example}

\begin{example}\label{eg:M12}
Another example is the sporadic group $M_{12}$ over $\bF_2$.
There are three simple modules $R$, $S$, and $T$ in the principal
block of $\bF_2M_{12}$, of dimensions
$44$, $1$ and $10$. We have $\dim_{\bF_2}\Ext^1_{\bF_2M_{12}}(R,T)=1$
and $\dim_{\bF_2}\Ext^1_{\bF_2M_{12}}(S,T)=2$, see Schneider~\cite{Schneider:1993a}.
\end{example}

\begin{example}
The following examples were computed using 
\textsf{Magma}~\cite{Bosma/Cannon/Playoust:1997a}.
For the group $PSL(3,\bF_4)$ over $\bF_4$, we can take for $R$, $S$
and $T$ simple modules of dimensions $8$, $9$ and $1$.
For the group $PSL(4,\bF_3)$ over $\bF_3$, we can take for $R$, $S$ and
$T$ simple modules of dimensions $19$, $1$ and $44$.
For the Higman--Sims group $HS$ over $\bF_2$, we can take for $R$, $S$
and $T$ the simple modules of dimensions $56$, $20$ and $1$.
\end{example}

\section{\texorpdfstring{One dimensional $\Ext^1$ groups}
{One dimensional Ext¹ groups}}\label{se:1Ext}

In some circumstances, our method can be modified to deal with groups
where all $\Ext^1$ groups between simple modules have dimension zero
or one. This technique works whenever there is an even length cycle in the
$\Ext^1$ quiver with alternating directions, and with some other arrow
head meeting two of its arrow heads:
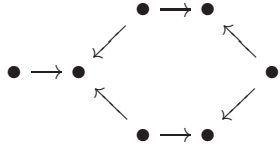
\[ \xymatrix@=4mm{&&\bullet\ar[r]\ar[dl]& \bullet\\
\bullet\ar[r]&\bullet&&&\bullet\ar[ul]\ar[dl]\\
&&\bullet\ar[ul]\ar[r]&\bullet}\] 
The method described in Section~\ref{se:2Ext} is a degenerate case of
this with a cycle of length two. This will lead to Case~(iii) of Theorem~\ref{th:Mm}.

Let $Q$ be the quiver above, with even cycle length $2r\ge 4$. 
Then a representation of $Q$ over $k$ consists of vector spaces $U$,
$V_1,\dots,V_r$ and $W_1,\dots,W_r$ together with $k$-linear maps
$X\colon U\to W_1$, $Y_i\colon V_i\to W_i$, 
$Z_{i+1}\colon V_{i+1}\to W_i$ and $Z_1\colon V_1\to W_r$. The case
$r=3$ is illustrated as follows.
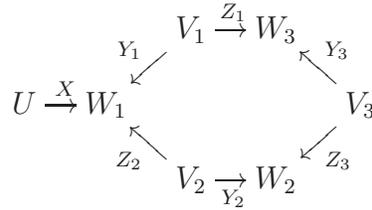
\vspace{-3mm}
\[ \xymatrix@=4mm{&&V_1\ar[r]^{Z_1}\ar[dl]_{Y_1}&W_3\\
U\ar[r]^X&W_1&&&V_3\ar[ul]_{Y_3}\ar[dl]^{Z_3}\\
&&V_2\ar[ul]^{Z_2}\ar[r]_{Y_2}&W_2}\]
The module $\sfM_m$ has basis vectors $v_1,\dots,v_m$ and
$w_1,\dots,w_m$. The space $U$ is spanned by $v_1$, $V_i$ is spanned
by the $v_j$ with $1<j\equiv i \pmod{r}$, and $W_i$ is spanned by the
$w_j$ with $1\le j\equiv i \pmod{r}$.
The actions of the arrows are as in
Table~\ref{eq:action}, with all $Y_i$ acting like $Y$ and all $Z_i$
acting like $Z$.  

\begin{theorem}\label{th:End1}
The endomorphism ring of the $kQ$-module $\sfM_m$ is equal to $k$, acting
as scalar multiples of the identity map.
\end{theorem}
\begin{proof}
The proof is the same as for Theorem~\ref{th:End3}.
\end{proof}

\begin{proof}[Proof of Theorem~\ref{th:Mm} in case {\rm (iii)}]
We have $kH$-modules $\Delta,\Delta_1,\dots,
\Delta_r$ with 
\[ \Rad(\Delta)\cong T_1,\ \Delta/\Rad(\Delta)\cong R,\ 
\Rad(\Delta_i)\cong T_i\oplus T_{i-1}\ (\text{or } T_1\oplus T_r
\text{ if } i=1),\ \Delta_i/\Rad(\Delta_i)\cong S_i. \] 
The $kH$-module $B$ is 
$\Delta\oplus\bigoplus_{i=1}^r\Delta_i\oplus\bigoplus_{i=1}^rT_i$, and
$\End_{kH}(B)\cong kQ$. This makes $B$ a
$kH$-$kQ$-bimodule with the property that as a right $kQ$-module it
is a projective generator for $\mmod$-$kQ$.
The functor $B \otimes_{kQ} -$ is fully faithful and exact. It sends
the simple $kQ$-modules to the corresponding simple $kH$-modules.
So this time, $M_m = B\otimes_{kQ}\sfM_m$ is a module of the following
form
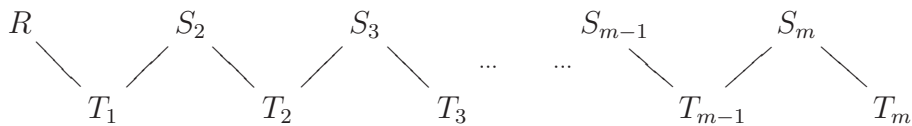
\[ \xymatrix@=5mm{R\ar@{-}[dr]&&
S_2\ar@{-}[dr]\ar@{-}[dl]&&
S_3\ar@{-}[dr]\ar@{-}[dl]&&\ar@{}[dl]|\cdots
&S_{m-1}\!\!\!\!\!\ar@{-}[dr]\ar@{}[dl]|\cdots&&
S_m\ar@{-}[dl]\ar@{-}[dr] \\
&T_1&&T_2&&T_3&&&T_{m-1}\!\!\!\!\!&&T_m} \]
where the indices are taken modulo $r$.
\end{proof}

\section{\texorpdfstring{Examples involving one dimensional $\Ext^1$ groups}
{Examples involving one dimensional Ext¹ groups}}\label{se:1Ext-eg}

In this section, we give examples satisfying Case~(iii) of Theorem~\ref{th:Mm}.

\begin{example}
The
$\Ext^1$ quiver of the Higman--Sims group $HS$ over a finite splitting field 
$k$ of characteristic three is as follows.
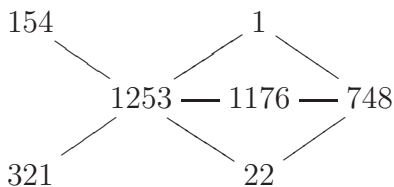
\[ \xymatrix@=5mm{154\ar@{-}[dr]&&1\ar@{-}[dr]\\
&1253\ar@{-}[ur]\ar@{-}[dr]\ar@{-}[r]&1176\ar@{-}[r]&748\\
321\ar@{-}[ur]&&22\ar@{-}[ur]} \]
This information comes from Section~2 of Waki~\cite{Waki:1993a}. The numbers
indicate the dimensions of the corresponding simple modules, and each
edge denotes a single arrow in each direction. 
Choosing the $4$-cycle $1\to 1253\leftarrow 22\to 748\leftarrow 1$ and
incoming arrow $154\to 1253$, the modules $M_m$ take the following form.
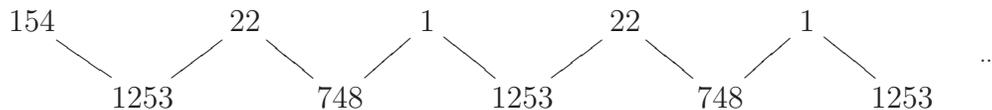
\[ \xymatrix@=5mm{
154\ar@{-}[dr]&&22\ar@{-}[dl]\ar@{-}[dr]&&
1\ar@{-}[dl]\ar@{-}[dr]&&
22\ar@{-}[dl]\ar@{-}[dr]&&
1\ar@{-}[dl]\ar@{-}[dr]&&\ar@{}[d]|\cdots\\
&1253&&748&&1253&&748&&1253&} \]
\end{example}

\begin{example}
Landrock and Michler~\cite{Landrock/Michler:1978a} showed that the
$\Ext^1$ quiver of the principal block of the sporadic Janko group $J_1$ over the splitting
field $\bF_4$ is as follows.
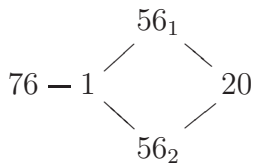
\[ \xymatrix@=3mm{&&56_1\ar@{-}[dl]\ar@{-}[dr]\\
76\ar@{-}[r]&1&&20\\
&&56_2\ar@{-}[ul]\ar@{-}[ur]} \]
So we can take the $4$-cycle $56_2\to 1 \leftarrow 56_1 \to
20\leftarrow 1$ and incoming arrow $76\to 1$.
\end{example}

\begin{example}
For $H=PSL(4,\bF_2)$ and $k=\bF_2$, the $\Ext^1$ quiver was
computed in Benson~\cite{Benson:1983a} to be as follows.
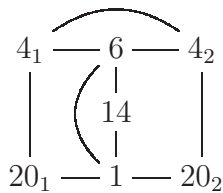
\[ \xymatrix@R=3mm@C=4mm{
4_1\ar@/^3ex/@{-}[rr]\ar@{-}[dd]\ar@{-}[r]
&6\ar@{-}[r]\ar@/_3ex/@{-}[dd]\ar@{-}[d]
&4_2\ar@{-}[dd]\\
&14&\\
20_1\ar@{-}[r]&1\ar@{-}[u]\ar@{-}[r]&20_2} \]
So for example we can take a $4$-cycle $1\to 6\leftarrow 4_1\to
20_1\leftarrow 1$
and incoming arrow $14 \to 6$. 
\end{example}

\begin{example}
Computations using \textsc{Magma}~\cite{Bosma/Cannon/Playoust:1997a} 
exhibit further examples. For the group
$H=PSL(3,\bF_5)$ and $k=\bF_5$, all the $\Ext^1$ groups are zero or
one dimensional, but the $\Ext^1$ quiver has a
$4$-cycle of the form $19\to 6 \leftarrow 18\to 8\leftarrow 19$ and an
incoming arrow $35 \to 6$.
\end{example}

\section{The construction}\label{se:constr}

In this section, we show how to pass from the constructions of
Sections~\ref{se:3Ext}, \ref{se:2Ext} and~\ref{se:1Ext} and groups satisfying the
hypotheses of Theorem~\ref{th:Mm} to examples proving
Theorem~\ref{th:main}. For convenience, we assume that $H$ is a finite
simple group, though this could probably be weakened somewhat.

Let $k$ and $Q$ be as in Section~\ref{se:3Ext},
\ref{se:2Ext}, or~\ref{se:1Ext}.
We define the $kQ$-module $\sfM_\infty$ to be the colimit of the
sequence of modules
\[ \sfM_2 \to \sfM_3\to\sfM_4 \to \cdots . \]
This has a basis consisting of the $v_i$ and $w_i$ with $i\ge 1$.
Then there are submodules $\sfM'_m$ of $\sfM_\infty$ such that 
$\sfM_m \cap \sfM'_m=0$, and
$\sfM_\infty/(\sfM_m + \sfM'_m)$ is one dimensional, spanned by the
image of $v_{i+1}$.
Namely, we take $\sfM'_m$ to be the submodule spanned by the
basis vectors $v_i$ with $i\ge m+2$ and $w_i$ with $i\ge m+1$.

Let $M_\infty$ and $M'_m$ be the $kH$-modules defined as the 
images of $\sfM_\infty$ and $\sfM'_m$
under $B \otimes_{kQ}-$. Thus $M_\infty$ is the colimit of the
sequence of modules
\[ M_2 \to M_3\to M_4 \to \cdots .\]
The submodules $M'_m$ of $M_\infty$ satisfy $M_m\cap M'_m=0$ and 
$M_\infty/(M_m+M'_m) \cong S$.

We define $G_m$ to be the group consisting of the 
vector space automorphisms $\alpha$ of $M_\infty$
such that $\alpha$ preserves and acts with determinant one on $M_m$,
and 
there exists $h \in H$ such that $\alpha$ acts in
the same way as $h$ on both $M_\infty/M_m$ and $M'_m$.
It follows from these conditions that $\alpha(M_i)\subseteq M_i$ 
for all $i\ge m$. 
Then since $M_m \subseteq M_{m+1}$ and 
$M'_{m+1}\subseteq M'_m$, we have $G_m\le G_{m+1}$.
The action of $H$ on $M_\infty$ gives an inclusion 
$i\colon H \to G_m$ for each $m$, compatible with these inclusions.

We have group homomorphisms $G_m \to H$ and $G_m \to SL(M_{m+1})$
given by the actions on $M'_m$ and on $M_{m+1}$.
Since $M'_m$ and $M_{m+1}$ span $M_\infty$, 
the intersection of the kernels of the actions on these spaces 
is trivial, so we obtain an injective homomorphism 
\[ G_m \to H \times SL(M_{m+1}). \]
This
is not surjective, because the image of $G_m$ in $SL(M_{m+1})$ 
stabilises $M_m$ setwise; but it properly contains $SL(M_m)$.
However, the image contains $H\times 1$, and is therefore the direct product
of $H$ and the image in $SL(M_{m+1})$. The subgroup mapping to
$H\times 1$ consists of elements acting as $h \in H$ on $M'_m$ and as
the identity of $M_\infty/M'_m$. This subgroup depends on $m$, 
and is not the same as the image of
$i\colon H\to G_m$. Indeed, the composite $H \xrightarrow{i} G_m \to
SL(M_{m+1})$ gives the action of $H$ on the module $M_{m+1}$.

\begin{lemma}\label{le:centraliser}
The centraliser of the image of $i\colon H \to G_m$ is trivial.
\end{lemma}
\begin{proof}
Let $z$ be an element of $C_{G_m}(i(H))$. Then $z$ commutes with the 
action of $H$ on $M_\infty$, and so $z$ acts as a scalar multiple of
the identity. But then $z$ also acts on $M_\infty/M_m$ as a multiple
of the identity. The action of $G_m$ on $M_\infty/M_m$ factors through
$G_m \to H$, and $Z(H)$ is trivial, so $z=1$.
\end{proof}

We define $G_\infty$ to be the colimit of the inclusions
\[ G_2 \to G_3 \to G_4 \to \cdots . \]
By Theorem~\ref{th:Mm}\,(3),
we have $C_{G_m}(H) \subseteq Z(G_m)=1$ and 
$C_{G_\infty}(H)\subseteq Z(G_\infty)=1$.\medskip

\begin{theorem}\label{th:semidirect}
The group $G_\infty$ is isomorphic to a semidirect product
$SL(\infty,k)\rtimes H$.
\end{theorem}
\begin{proof}
By the definition of $G_m$, there is a homomorphism $G_m\to H$, taking
an element $\alpha$ to the element $h\in H$ that acts on $M_\infty/M_m$ and
$M'_m$ in the same way as $\alpha$. The composite $H \xrightarrow{i} G_m \to H$ is
the identity.
These homomorphisms are compatible
with the inclusions $G_m\to G_{m+1}$, and therefore describe a
homomorphism $G_\infty \to H$ with kernel $SL(\infty,k)$ a normal
complement to the inclusion $i\colon H\to G_\infty$.
\end{proof}

\begin{theorem}\label{th:action}
Given a finite simple group $H$ satisfying condition {\rm (i)}, {\rm (ii)} or
{\rm (iii)} of
Theorem~\ref{th:Mm}, there is an action of $H$ on $SL(\infty,k)$
with trivial centraliser.
\end{theorem}
\begin{proof}
This follows from Theorem~\ref{th:semidirect} and
Lemma~\ref{le:centraliser}.
\end{proof}

\begin{proof}[Proof of Theorem~\ref{th:main}]
We set $k=\bF_q$ in Theorem~\ref{th:action}.
For $p$ odd, we can use Example~\ref{eg:p2}, and for $p=2$ we can use
Example~\ref{eg:M12} for the choice of $H$; there are, of course, many
other possible choices.
\end{proof}

\bibliographystyle{amsplain}
\bibliography{../repcoh}

\end{document}